\documentclass[12pt]{article}
\begin{document}
%\global\baselineskip24pt plus2pt minus2pt

\catcode`ð=\active
\defð{\u{g}}
\catcode`Ð=\active
\defÐ{\u{G}}
\catcode`Ý=\active
\defÝ{\. I}
\catcode`ö=\active
\defö{\"{o}}
\catcode`Ö=\active
\defÖ{\"O}
\catcode`ü=\active
\defü{\"{u}}
\catcode`Ü=\active
\defÜ{\"{U}}
\catcode`Þ=\active
\defÞ{\c{S}}
\catcode`þ=\active
\defþ{\c{s}}
\catcode`ý=\active
\defý{{\i}}
\catcode`ç=\active
\defç{\d{c}}
\catcode`Ç=\active
\defÇ{\d{C}}
\def\oast{\mathbin{\bigcirc\llap{$\ast$\kern.1cm}}}
\def\Nlim{\mathop{\rm N\!-\!lim}}
\def\nLOP{\mathop{\rm lim}}
\def\nn{\nLOP\limits_{n\to\infty}}
\def\ne{\nLOP\limits_{\varepsilon\to 0}}
\def\NLOP{\mathop{\rm N\!-\!lim}}
\def\Nn{\NLOP\limits_{n\to\infty}}
\def\Nm{\NLOP\limits_{m\to\infty}}
\def\Nv{\NLOP\limits_{\nu\to\infty}}
\def\Ne{\NLOP\limits_{\varepsilon\to 0}}
\def\cosec{\mathop{\rm cosec}\nolimits}
\def\coshec{\mathop{\rm coshec}\nolimits} \def\Ci{\mathop{\rm Ci}\nolimits}
\def\Si{\mathop{\rm Si}\nolimits} \def\oot{\mathop{\textstyle {1\over
2}}\nolimits} \def\fr#1#2{{\textstyle {#1\over#2}}}
\def\sgn{\mathop{\rm sgn}\nolimits} \def\itr{\mathop{\rm
int}\nolimits} \def \tast{\mathbin{\bigtriangleup \!\!\!\!*\,}}
\def\bast{\mathbin{\vcenter{\hrule\hbox to 8.3pt{\vrule\kern0.9pt\vbox{\hbox
{\mathstrut$\ast$}}\kern0.9pt\vrule}\hrule}}}
\def\lamb{\mathop{\lambda\!\!\!\!\lambda}}
\newcommand{\bsym}[1]{\mbox{\boldmath$#1$}}
\def\mub{\mathop{\mu\!\!\!\!\mu}}
\def\muf{\mathop{\mu\!\!\!\!\mu}}
\newcommand{\ol}{\overline}
\newcommand{\ul}{\underline}
\newcommand{\bc}{\begin{center}}
\newcommand{\ec}{\end{center}}
\newcommand{\ts}{\textstyle}
\newcommand{\ds}{\displaystyle}
\newcommand{\be}{\begin{equation}}
\newcommand{\ee}{\end{equation}}
\newcommand{\beq}{\begin{eqnarray}}
\newcommand{\eeq}{\end{eqnarray}}
\newcommand{\beqa}{\begin{eqnarray*}}
\newcommand{\eeqa}{\end{eqnarray*}}
\newcommand{\lam}{\lambda}
\def\lamb{\mathop{\lambda\!\!\!\!\lambda}}
\newcommand{\vphi}{\varphi}
\newcommand{\Aa}{{\cal A}}
\newcommand{\BB}{{\cal B}}
\newcommand{\DD}{{\cal D}}
\newcommand{\EE}{{\cal E}}
\newcommand{\GG}{{\cal G}}
\newcommand{\RR}{{\cal R}}
\newcommand{\TT}{{\cal T}}
\newcommand{\XX}{{\cal X}}
\newcommand{\YY}{{\cal Y}}
\newcommand{\ZZ}{{\cal Z}}
\newcommand{\CC}{{\cal C}}
\newcommand{\n}{\nu}
\newcommand{\la}{\langle}
\newcommand{\ra}{\rangle}

\medskip
\bc \large {{\bf  Applications of Neutrix Calculus to Special
Functions in Conjunction with Polygamma Functions }} \ec

\bc {\bf Emin Özçað} \ec

\bigskip \noindent
{\bf Abstract.} {\footnotesize In this paper we define the
polygamma functions $\psi^{(n)}(x)$ for negative integers by using
neutrix calculus. }



\bigskip \bigskip \noindent
\bc {\large {\bf 1. Introduction}} \ec

\bigskip

\medskip \noindent
A rich combination of elementary and special functions arises in
the evaluation of Euler sums. Special functions which typically
appear are the gamma, beta, digamma and polygamma functions, the
zeta functions, polylogarithms, hypergeometric functions, and
logarithmic-trigonometric integrals.

\medskip \noindent
Kirchoff was first to apply the polygamma functions in physics,
the summation of rational series and the evaluation of integrals
are some applications that are still relevant. Recently, the
summation of series containing $\psi^{(n)}(z)$ was arisen in
Feynman calculations \cite{kn:coff}.

\medskip \noindent
Kölbig presented two formulae for the $\psi^{(n)}(p/q)$ by
using the series definition of polylogarithm function, see
\cite{kn:kol1,kn:kol2}. Coffey considered the sums over the
digamma function, containing summand with $(\pm
1)^n\psi(n+p/q)/n^2$ and made extension to sums over the polygamma
functions in \cite{kn:coff}.

\bigskip \noindent The technique of neglecting appropriately
defined infinite quantities was devised by Hadamard and the
resulting finite value extracted from the divergent integral is
usually referred to as the Hadamard finite part, see
\cite{kn:jones}.

\noindent Using the concepts of the neutrix and the neutrix limit,
Fisher gave the general principle for the discarding of unwanted
infinite quantities from asymptotic expansions and has been
exploited in context of distributions and special functions, see
\cite{kn:f1,kn:f2,kn:fbil,kn:fkurb1,kn:ozinci1,kn:ozinci2}

\medskip \noindent
Y. Jack Ng and H. van Dam applied the neutrix calculus, in
conjuction with the Hadamard integral, developed by van der Corput
see \cite{kn:c}, to quantum field theories, in particular, to
obtain finite results for the cofficients in the perturbation
series. They also applied neutrix calculus to quantum field
theory, obtaining finite renormalization in the loop calculations,
see \cite{kn:Ng1,kn:Ng2}

\medskip \noindent
In the following we let $N$ be the neutrix \cite{kn:c} having
domain $N' = {\{\varepsilon:\,0<\varepsilon<\infty}\}$ and range
$N''$ the real numbers, with negligible functions finite linear
sums of the functions \be \varepsilon^{\lambda} \ln^{r-1}
\varepsilon, \quad \ln^r \varepsilon \qquad (\lambda <0, \quad  r
= 1,2,\ldots) \ee and all functions $f(\varepsilon)$ which
converge to zero in the normal sense as $\varepsilon$ tends to
zero.

\medskip \noindent
If $f(\varepsilon)$ is a real (or complex) valued function defined
on $N'$ and if it is possible to find a constant $c$ such that
$f(\varepsilon) -c$ is in $N$, then $c$ is called the neutrix
limit of $f(\varepsilon)$ as $\varepsilon \rightarrow 0$ and we
write $ \Nlim_{\varepsilon\to 0} f(\varepsilon) =c.$

\medskip \noindent
Note that if a function $f( \varepsilon)$ tends to $c$ in the
normal sense as $\varepsilon$ tends to zero, it converges to $c$
in the neutrix sense.

\medskip \noindent
Also note that if a function $H(\varepsilon)=\upsilon(\varepsilon
)+f(\varepsilon ),$ where $\upsilon(\varepsilon )$ is the sum of
negligible functions of $H(\varepsilon ),$ then p.f.
$H(\varepsilon ),$ Hadamard's finite part of $H(\varepsilon ),$ is
equal to $f(\varepsilon )$ and so $$\Ne H(\varepsilon )=\ne
f(\varepsilon )=\ne  p.f. H(\varepsilon ).$$

\noindent The reader may find the general definition of the
neutrix limit with some examples in \cite{kn:c,kn:f1,kn:f2}.

\medskip \noindent
In this paper we use Fisher's principle to define the derivative
of the digamma function for negative integers. First of all, we
give the definition of the gamma function for all $x.$

\medskip \noindent
The gamma function $\Gamma (x)$ is usually defined for $x>0$ by
\beq \Gamma(x) =\int_0^\infty t^{x - 1}e^{-t} \,dt \eeq the
integral only converging for $x>0,$ see \cite{kn:gs,kn:rr}. It
follows from equation (2) that
\be \Gamma (x +1) =x\Gamma(x) \ee
for $x>0$ and this equation is used to define $\Gamma (x)$ for
negative, non-integer values of x. Using the regularization,
Gelfand and Shilov \cite{kn:gs} define the gamma function
$$\Gamma (x) =\int_0^1 t^{x-1}\Bigl [e^{-t} -\sum_{i=0}^{n-1}(-1)^i{t^i\over
i!}\Bigr ]dt +\int_1^\infty t^{x-1} e^{-t} \,dt +\sum_{i=0}^{n-1}
{(-1)^i\over i! (x+i)} $$ for $x>-n, \quad x\neq
0,-1,-2,\ldots,-n+1$ and
$$\Gamma (x) = \int_0^\infty t^{x-1}\Bigl [ e^{-t} -\sum_{i=0}^{n-1}(-1)^i
{t^i\over i!}\Bigr ]\,dt  $$ for $-n < x <-n+1.$

\medskip \noindent
Fisher proved that
$$ \Gamma (x) = \Ne
\int_{\varepsilon}^{\infty}t^{x - 1}\,e^{-t} \,dt
$$ $x \neq 0,-1,-2,\ldots$ and defined $\Gamma(-m)$ by
\beq &&\Gamma (-m) =
\Ne \int_{\varepsilon}^{\infty}t^{-m- 1}\,e^{-t}
\,dt \hskip2.5in \nonumber \\
&&\hskip0.2in = \int_1^\infty t^{-m - 1}e^{-t}\,dt +\int_0^1 t^{-
m -1}\Bigl [ e^{-t} - \sum_{i=0}^m\frac{(-t)^i}{i!} \Bigr ]\,dt
-\sum_{i=0}^{m-1} \frac{(-1)^i}{i! (m-i)} \eeq for $m=1,2,\ldots,$
see \cite{kn:f4}.

\medskip \noindent More generally, the r-th
derivative of gamma function $\Gamma (x)$ is defined by \be
\Gamma^{(r)} (x) = \Ne \int_{\varepsilon}^{\infty}t^{x - 1} \ln
^rt \,e^{-t} \,dt \ee for all $x$ and $r=0,1,2,\ldots, $ see
\cite{kn:fkurb1}.

\medskip \noindent
Fisher obtained from equation (4) that
$$ \Gamma(-m)+\frac{1}{m}\Gamma(-m+1)=\frac{(-1)^m}{mm!} $$
from which it followed by induction that \be
\Gamma(-m)=\frac{(-1)^m}{m!}[ \phi(m)-\gamma ] \ee for
$m=1,2,\ldots, $ where $\gamma $ denotes Euler's constant and
$$\phi(m)=\left\{ \begin{array}{cc} 0,& m=0,\\ \sum_{i=1}^m
\frac{1}{i},& m=1,2,\ldots \end{array}\right.$$ in particular
$\Gamma(0)=\Gamma'(1)=-\gamma ,$ see \cite{kn:f4}.

\bigskip \noindent The digamma function $\psi(x)=\frac{d}{dx}\log\Gamma(x)=
\frac{\Gamma'(x)}{\Gamma(x)}$ defined for $x>0$ has the integral
representation
$$\psi(x)=-\gamma+\int_0^\infty \frac{e^{-t}-e^{-xt}}{1-e^{-t}}\,dt $$
and it can be written as \be \psi(x)= -\gamma+\int_0^1
\frac{1-t^{x-1}}{1-t}\,dt \qquad (x>0). \ee

\medskip \noindent
Differentiating equation (3), we have \be
\Gamma'(x+1)=\Gamma(x)+x\Gamma'(x) \ee for $x\neq 0,-1,-2,\ldots$
and it follows that
$$\psi(x+1)=\psi(x)+\frac{1}{x} $$
see \cite{kn:ln,kn:moll} and this equation is used to define
$\psi(x)$ for negative, noninteger values of $x.$ Thus if
$-m<x<-m+1; \quad m=1,2,\ldots,$ then \be \psi(x)=
-\gamma+\int_0^1
\frac{1-t^{x-1+m}}{1-t}\,dt-\sum_{k=1}^{m-1}\frac{1}{x+k} \ee

\medskip \noindent
In \cite{kn:fkurb2} Fisher and kuribayashi proved the following
equations to define $\psi(-m)$ for $m=1,2,\ldots. $

\medskip \noindent {\bf Theorem 1.1.}
\beqa \Gamma^{(n)}(x)&=& \Ne \Gamma^{(n)}(x+\varepsilon)\\
\Gamma(x+1)&=& \Ne (x+\varepsilon)\Gamma(x+\varepsilon) \\
\Gamma'(x+1) &=& \Ne \bigl
[\Gamma(x+\varepsilon)+(x+\varepsilon)\Gamma'(x+\varepsilon)\bigr
] \eeqa for all $x$ and $n=0,1,2,\ldots.$

\medskip \noindent
Theorem 1.1 suggested that the digamma function $\psi(-m)$ can be
defined by
$$\psi(-m)=\Ne \frac{\Gamma'(-m+\varepsilon)}{\Gamma(-m+\varepsilon)}$$
for $m=0,1,2,\ldots,$ provided the neutrix limit exists, and with
this definition

\be \psi(-m)=\psi(1)+\phi(m)=-\gamma+\phi(m) \ee for
$m=0,1,2,\dots,$ see \cite{kn:fkurb2}.

\medskip \noindent
Recently Tuneska and Jolevski used the integral representation of
the digamma function to obtain same result given in
\cite{kn:biljana}.

\bigskip\noindent \bc {\large {\bf 2. Defining Polygamma  Function
$\psi^{(n)}(-m)$ }} \ec

\bigskip \noindent
The polygamma function is defined by \be
\psi^{(n)}(x)=\frac{d^n}{dx^n}\psi(x)=\frac{d^{n+1}}{dx^{n+1}}\ln
\Gamma(x) \qquad (x>0). \ee
It may be represented as $$
\psi^{(n)}(x)= (-1)^{n+1}\int_0^\infty \frac{t^n
e^{-xt}}{1-e^{-t}}$$ which holds for $x>0,$ and \be \psi^{(n)}(x)=
-\int_0^1\frac{t^{x-1}\ln ^nt}{1-t}\,dt. \ee

\medskip \noindent
It satisfies the recurrence relation \be
\psi^{(n)}(x+1)=\psi^{(n)}(x)+\frac{(-1)^nn!}{x^{n+1}} \ee see
\cite{kn:btr,kn:kol1,kn:kol2,kn:ln,kn:rr}. This is used to define
the polygamma function for negative non-integer values of $x.$
Thus if $-m<x<-m+1, \quad m=1,2,\ldots,$ then \be
\psi^{(n)}(x)=-\int_0^1 \frac{t^{x+m-1}}{1-t}\ln^n
t\,dt-\sum_{k=0}^{m-1}\frac{(-1)^nn!}{(x+k)^{n+1}}.\ee

\medskip \noindent
Kölbig gave the formulae for the integral $\int_0^1 t^{\lambda
-1}(1-t)^{-\nu}\ln ^mt\,dt$ for integer and half-integer values of
$\lambda $ and $\nu $ in \cite{kn:kol3}. As the integral
representation of the polygamma function is similar to the
integral mentioned above, by using the neutrix limit we prove the
existence of the integral in (12) as follows.

\medskip \noindent
Now we let $N$ be a neutrix having domain the open interval
$\{ \epsilon : 0<\varepsilon < {1\over 2}\}$ with the same
negligible functions as in equation (1). We first of all need the
following lemma.

\medskip \noindent
{\bf Lemma 2.1} The neutrix limits as $\varepsilon $ tends to zero
of the functions $$ \int_{\varepsilon}^{1/2} t^x \ln ^nt\ln
^r(1-t)\,dt, \qquad \int_{1/2}^{1-\varepsilon}(1-t)^x\ln ^nt\ln
^r(1-t)\,dt $$ exists for $n,r=0,1,2,\ldots$ and all $x.$

\medskip \noindent
{\bf Proof.} Suppose first of all that $n=r=0.$ Then
$$ \int_\varepsilon ^{1/2} t^x \,dt = \left \{ \begin{array}{cc}
\frac{2^{-x-1}-\varepsilon^{x+1}}{x+1}, & x\neq -1, \\
-\ln 2-\ln \varepsilon, & x=-1
\end{array} \right. $$
and so $\Ne \int_\varepsilon^{1/2} t^x \,dt $ exists for all $x.$

\noindent Now suppose that $r=0$ and that $ \Ne
\int_\varepsilon^{1/2} t^x\ln ^n t \,dt $ exists for some
nonnegative integer n and all $x.$
Then
$$ \int_\varepsilon ^{1/2} t^x\ln^{n+1} t \,dt = \left \{ \begin{array}{cc}
\frac{-2^{-x-1}\ln^{n+1} 2-\varepsilon ^{x+1}\ln^{n+1}
\varepsilon}{x+1}- \frac{n+1}{x+1}\int_\varepsilon ^{1/2} t^x\ln^n
t \,dt, & x\neq -1,
\\ \\ \frac{(-1)^n\ln^{n+2} 2-\ln^{n+2} \varepsilon }{m+2}, & x=-1
\end{array} \right.$$
and it follows by induction that $ \Ne \int_\epsilon^{1/2} t^x
\ln^n t \,dt $ exists for $n=0,1,2,\ldots$ and all $x.$

\noindent Finally we note that we can write
$$ \ln^r (1-t) = \sum_{i=1}^\infty \alpha_{in} t^i $$
for $r=1,2,\ldots,$ the expansion being valid for $|t|<1.$
Choosing a positive integer k such that $x+k>-1,$ we have \beqa
&&\int_\varepsilon ^{1/2} t^x \ln^n t \ln^r(1-t) \,dt =\hskip2.8in \\
&&\hskip0.2in = \sum_{i=1}^{k-1}\alpha_{in}\int_\varepsilon ^{1/2}
t^{x+i} \ln^n t\,dt +\sum_{i=k}^\infty
\alpha_{in}\int_\varepsilon^{1/2} t^{x+i} \ln^n t\,dt. \eeqa

It follows from what we have just proved that
$$\Ne \sum_{i=1}^{k-1}\alpha_{in}\int_\varepsilon ^{1/2} t^{x+i}
\ln^n t\,dt $$ exists and further \beqa \Ne \sum_{i=k}^\infty
\alpha_{in}\int_\varepsilon ^{1/2} t^{x+i} \ln^n t\,dt &=& \ne
\sum_{i=k}^\infty \alpha_{in}
\int_\varepsilon ^{1/2} t^{x+i} \ln^n t\,dt \\
&=& \sum_{i=k}^\infty \alpha_{in} \int_0^{1/2} t^{x+i} \ln^n
t\,dt, \eeqa proving that the neutrix limit of $\int_\varepsilon
^{1/2} t^x \ln^n t \ln^r(1-t) \,dt$ exists for $n,r=0,1,2,\ldots$
and all $x.$ Making the substitution $1-t=u$ in
$$\int_{1/2}^{1-\varepsilon }
(1-t)^x \ln ^n t \ln ^r t\,dt,$$ it follows that
$\int_{1/2}^{1-\varepsilon } (1-t)^x \ln ^n t \ln ^r t\,dt$ also
exits for $n,r=0,1,2,\ldots$ and all $x.$

\medskip \noindent
{\bf Lemma 2.2} The neutrix limit as $\varepsilon \rightarrow 0$
of the integral $\int_{\varepsilon }^1 t^{-m-1} \ln ^n t\,dt$
exists for $m,n=1,2,\ldots$ and \be \Ne \int_{\varepsilon }^1
t^{-m-1} \ln ^n t\,dt =-\frac{n!}{m^{n+1}}. \ee

\medskip \noindent
{\bf Proof.} Integrating by parts, we have
$$\int_{\varepsilon}^1 t^{-m-1} \ln t\,dt = m^{-1}\varepsilon^{-m} \ln \varepsilon
+m^{-1}\int_{\varepsilon}^1 t^{-m-1} \,dt $$ and so $$\Ne
\int_{\varepsilon}^1 t^{-m-1} \ln t\,dt = -\frac{1}{n^2} $$
proving equation (15) for $n=1$ and $m=1,2,\ldots.$ Now assume
that equation (15) holds for some $m$ and $n=1,2,\ldots.$ Then
\beqa \int_{\varepsilon}^1 t^{-m-2} \ln ^n t\,dt &=&
(m+1)^{-1}\varepsilon^{-m-1}\ln ^n\varepsilon
+\frac{n}{m+1}\int_{\varepsilon}^1 t^{-m-2} \ln ^{n-1} t\,dt \\
&=& (m+1)^{-1}\varepsilon^{-m-1}\ln ^n\varepsilon +
\frac{n}{m+1}\frac{-(n-1)!}{(m+1)^n} \eeqa and it follows that
$$\Ne \int_{\varepsilon}^1 t^{-m-2} \ln ^n t\,dt
=-\frac{n!}{(m+1)^{n+1}}$$ proving equation (15) for $m+1$ and
$n=1,2,\ldots.$

\medskip\noindent
Using the regularization and the neutrix limit, we prove
the following theorem.

\medskip \noindent
{\bf Theorem 2.3} The function $\psi^{(n)}(x)$ exists for
$n=0,1,2,\ldots,$ and all $x.$

\medskip \noindent
{\bf Proof.} Choose positive integer r such that $x>-r.$ Then we
can write \beqa \int_\varepsilon ^{1-\varepsilon
}\frac{t^{x-1}}{1-t}\ln^n t\,dt &=& \int_\varepsilon ^{1/2}
t^{x-1}\ln^n t \Bigl [\frac{1}{1-t}-\sum_{i=0}^{r-1}(-1)^i t^i
\Bigr ]\,dt \\
&& \hskip0.1in + \sum_{i=0}^{r-1}(-1)^i \int_\varepsilon ^{1/2}
t^{x+i-1}\ln^n t\,dt+\int_{1/2}^{1-\varepsilon }
\frac{t^{x-1}}{1-t}\ln^n t\,dt. \eeqa

We have \beqa && \ne \int_\varepsilon
^{1/2}\frac{t^{x-1}}{1-t}\ln^n t\Bigl
[\frac{1}{1-t}-\sum_{i=0}^{r-1}(-1)^i t^i \Bigr ]\,dt =
\hskip2.0in
\\ &&\hskip1.7in =\int_0^{1/2} \frac{t^{x-1}}{1-t}\ln^n t \Bigl
[\frac{1}{1-t}-\sum_{i=0}^{r-1}(-1)^i t^i \Bigr ]\,dt \eeqa and
$$\ne \int_{1/2}^{1-\varepsilon }
\frac{t^{x-1}}{1-t}\ln^n t\,dt = \int_{1/2}^1
\frac{t^{x-1}}{1-t}\ln^n t\,dt$$ the integrals being convergent.
Further, from the Lemma 2.1 we see that the neutrix limit of the
function $$\sum_{i=0}^{r-1}(-1)^i \int_\varepsilon ^{1/2}
t^{x+i-1}\ln^n t\,dt $$ exists and implying that
$$\Ne \int_\varepsilon ^{1-\varepsilon
}\frac{t^{x-1}}{1-t}\ln^n t\,dt $$ exists. This proves the
existence of the function $\psi^{(n)}(x)$ for $n=0,1,2,\ldots,$
and all $x.$

\medskip \noindent
Before giving our main theorem, we note that
$$\psi^{(n)}(x)=-\Ne \int_\varepsilon
^1\frac{t^{x-1}}{1-t}\ln^n t\,dt$$ since the integral is
convergent in the neighborhood of the point $t=1.$

\medskip \noindent
{\bf Theorem 2.4} The function $\psi^{(n)}(-m)$ exists and \be
\psi^{(n)}(-m)= \sum_{i=1}^m\frac{n!}{i^{n+1}}+(-1)^{n+1}n!\zeta
(n+1)\ee for $n=1,2,\ldots$ and $m=0,1,2,\ldots,$ where $\zeta(n)$
denotes zeta function.

\medskip \noindent
{\bf Proof.} From Theorem 2.3, we have \beq \psi^{(n)}(-m)&=& -\Ne
\int_\varepsilon^1\frac{t^{-m-1}\ln ^n t}{1-t}\,dt \nonumber \\
&=&-\Ne \int_\varepsilon^1\Bigl [\sum_{i=1}^{m+1}
t^{-i}+(1-t)^{-1}\Bigr ]\ln ^nt\,dt. \eeq We first of all evaluate
the neutrix limit of integral $\int_\varepsilon^1 t^{-i}\ln^n
t\,dt$ for $i=1,2,\ldots$ and $n=1,2,\ldots.$

\noindent
It follows from Lemma 2.2 that \be \Ne
\int_\varepsilon^1 t^{-i}\ln^n t\,dt=-\frac{n!}{(i-1)^{n+1}}\qquad
(i>1).\ee

\noindent For $i=1,$ we have
$$\int_\varepsilon^1 t^{-1}\ln^n t\,dt = O(\varepsilon ).$$
Next \beq \int_\varepsilon^1 \frac{\ln ^n t}{1-t}\,dt &=&\int_0^1
\frac{\ln ^n t}{1-t}\,dt=\sum_{k=0}^\infty \int_0^1 t^k\ln ^nt\,dt
\nonumber \\ &=&\sum_{k=0}^\infty \frac{(-1)^n n!}{(k+1)^{n+1}}
\nonumber \\ &=& (-1)^n n!\zeta(n+1), \eeq where $$\int_0^1 t^k\ln
^nt\,dt = \frac{(-1)^n n!}{(k+1)^{n+1}}.$$ It now follows from
equations (17),(18) and (19) that \beqa \psi^{(n)}(-m)&=& -\Ne
\int_\varepsilon^1\frac{t^{-m-1}\ln ^n t}{1-t}\,dt  \\
&=& -\sum_{i=1}^{m+1}\Ne \int_\varepsilon^1 t^{-i}\ln^n t\,dt-\Ne
\int_\varepsilon^1 \frac{\ln^n t}{(1-t)}\,dt. \\
&=& \sum_{i=1}^m\frac{n!}{i^{n+1}}+(-1)^{n+1}n!\zeta (n+1) \eeqa
implying equation (16).

\bigskip \noindent  Note that the digamma function $\psi(x)$ can be
defined by \be \psi (x)= -\gamma+\Ne \int_{\varepsilon}^1
\frac{1-t^{x-1}}{1-t}\,dt \ee for all $x.$

\medskip \noindent
Using Lemma 2.1 we have
\beqa \int_{\varepsilon}^1
\frac{1-t^{-m-1}}{1-t}\,dt
&=&-\sum_{i=1}^{m+1}\int_{\varepsilon}^1
t^{-i}\,dt \\
&=&-[\ln 1-\ln
\varepsilon]-\sum_{i=2}^{m+1}\frac{[1-\varepsilon^{-i+1}]}{-i+1}
\eeqa and it follows from equation (20) that
\beqa \psi (-m) &=&
-\gamma+\Ne \int_{\varepsilon}^1
\frac{1-t^{-m-1}}{1-t}\,dt=-\gamma
+\sum_{i=1}^m i^{-1} \\
&=&-\gamma +\phi (m) \eeqa which was obtained in \cite{kn:fkurb2}
and \cite{kn:biljana}.

\bigskip \noindent
Emin Özçað \\
Department of mathematics \\
Hacettepe University Beytepe 06800 \\
Ankara, Turkey \\
e-mail : ozcag1@hacettepe.edu.tr


\bigskip
\begin{thebibliography}{99}

\bigskip
\bibitem{kn:btr} Batýr, N., On some properties of digamma and
polygamma functions, J. Math. Anal. Appl., 328(2007) 452-465.
\bibitem{kn:coff} Coffey, M. W., On one-dimensional
digamma and polygamma series related to the evaluation of Feynman
diagrams, J. Comput. Appl. Math., 183(2005) 84-100.
\bibitem{kn:c} van der Corput, J. G., Introduction to the neutrix
calculus, J. Analyse Math., 7(1959) 291-398.
%\bibitem{kn:davi} Davis, H. T., An extension to polygamma
%functions of a theorem of Gauss, Bull. Am. Math. Soc., 41(1935)
%243-247.
\bibitem{kn:f1} Fisher, B., Neutrices and the product of distributions,
Studia Math., 57(1976) 263-274.
\bibitem{kn:f2} Fisher, B., A non-commutative neutrix product of
distributions, Math. Nachr., 108(1982) 117-127.
%\bibitem{kn:f3} Fisher, B., Neutrices  and distributions, Fourt Internat. Conf.
%on Complex Analysis and Applications, Varna, Bulgaria, (1987),
%169-175.
\bibitem{kn:f4} Fisher, B., On defining $\Gamma(-n)$ for
$n=0,1,2,\ldots $, Rostock. Math. Kolloq., 31(1987) 4-10.
%\bibitem{kn:f5} Fisher, B., On defining the incomplete gamma function $\gamma (-m,x_-),$
%Integral Trans. Spec. Funct., 15(6)(2004) 467-476.
\bibitem{kn:fbil} Fisher, B., Jolevska-Tuneska, B., and Kýlýçman, A.,
On defining the incomplete Gamma function, Integral Trans. Spec.
Funct., 14(4)(2003) 293-299.
\bibitem{kn:fkurb1} Fisher, B. and Kuribayashi, Y., Neutrices and the
Gamma function, J. Fac. Ed. Tottori Univ. Mat. Sci., 36(1-2)(1987)
1-7.
\bibitem{kn:fkurb2} Fisher, B. and Kuribayashi, Y., Some results on the
Gamma function, J. Fac. Ed. Tottori Univ. Mat. Sci., 3(2)(1988)
111-117.
\bibitem{kn:gs}I. M. Gel'fand and G. E. Shilov,
Generalized Functions, Vol. I., Academic Press, Newyork/London,
(1964).
\bibitem{kn:Ng1}Ng, Jack Y. and van Dam, H., Neutrix Calculus and
Finite Quantum Field Theory, J. Phys., A; Math. Gen. 38(2005)
317-323.
\bibitem{kn:Ng2} Ng, Jack Y. and van Dam, H., An application of neutrix
calculus to Quantum Field Theory, Internat. J. Modern Phys.,
A21(2)(2006) 297-312.
\bibitem{kn:biljana} Jolevska-Tuneska, B. and Jolevski, I., Some
results on the digamma function, Appl. Math. Inf. Sci., 7(2013)
167-170.
\bibitem{kn:jones} Jones, D. S., Hadamard's Finite Part, Math.
Methods Appl. Sci., 19(1996) 1017-1052.
\bibitem{kn:kol3} Kölbig, K. S., On the integral
$\int_0^1 x^{\nu -1}(1-x)^{-\lambda}\ln ^mx\,dx,$ J. Comput. Appl.
Math., 18(1987) 369-394.
\bibitem{kn:kol1} Kölbig, K. S., The polygamma function and the
derivatives of the cotangent function for rational arguments,
CERN-IT-Reports, CERN-CN-96-005, 1996.
\bibitem{kn:kol2} Kölbig, K. S.,The polygamma function
$\psi ^{(k)}(x)$ for $x=1/4$ and $x=3/4$, J. Comput. Appl. Math.,
75(1996) 43-46.
\bibitem{kn:ln} Laforgia, A. and Natalini, P., Exponentials, gamma
and polygamma functions: Simple proofs of classical and new
inequalities, J. Math. Anal. Appl. 407(2013) 497-504.
\bibitem{kn:moll} Medina, A. L. and Moll, V. H., The integrals in
Gradshteyn and Ryzhik,. Part 10: The digamma function, Scientia,
Series A: Math. Sci., 17(2009) 45-66.
\bibitem{kn:ozinci1} Özçað, E., Ege, Ý., Gürçay, H. and
Jolevska-Tuneska, B., Some remarks on the incomplete gamma
function, in: Kenan Ta\c s et al. (Eds.), Mathematical Methods in
Engineering, Springer, Dordrecht, 2007, pp. 97-108.
\bibitem{kn:ozinci2} Özçað, E., Ege, Ý., Gürçay, H. and
Jolevska-Tuneska, B., On partial derivatives of the incomplete
beta function, Appl. Math. Lett., 21(2008) 675-681.
\bibitem{kn:rr} Rradshhteyn, I. S. and Ryzhik, I. M., Tables of integrals, Series, and
Products, Academic Press, San Diego, 2000.

\end{thebibliography}
\end{document}